\documentclass{elsarticle}

\usepackage{mathtools}
\usepackage{float}
\usepackage[utf8]{inputenc}   
\usepackage[T1]{fontenc}
\usepackage{amsthm}

\usepackage{natbib}

\newtheorem{thm}{Theorem}
\newtheorem{cor}{Corollary}
\newtheorem{prop}{Proposition}

\newtheorem{definition}{Definition}

\usepackage{hyperref}
\usepackage{dsfont}
\usepackage{lipsum}
\usepackage{mathtools}
\usepackage{amssymb}
\usepackage{float}
\usepackage{caption}
\usepackage{geometry}

\newcommand{\D}{\mathcal{D}} 

\newcommand{\type}{\opname{type}}

\newcommand{\opname}[1]{\operatorname{\textup{\textsf{#1}}}}

\newcommand{\oplusI}{\oplus_{\mathcal{I}}}
\newcommand{\oplusF}{\oplus_{\mathcal{F}}}
\newcommand{\otimesI}{\otimes_{\mathcal{I}}}
\newcommand{\otimesF}{\otimes_{\mathcal{F}}}

\newcommand{\size}{\opname{size}}

\newcommand{\Area}{\opname{area}}

\newcommand{\Conj}{\opname{Conj}}
\newcommand{\ConjD}{\Conj_{\D}}
\newcommand{\ConjI}{\Conj_{\mathcal{I}}}

\newcommand{\ConjF}{\Conj_{\mathcal{F}}}

\newcommand{\lat}{\opname{lat}}
\newcommand{\Long}{\opname{long}}

\newcommand{\FF}{\mathrm}

\newcommand{\ljaw}{\opname{jaw}}

\newcommand{\Val}{\opname{Val}}
\newcommand{\Peak}{\opname{Peak}}
\newcommand{\DR}{\opname{DR}}
\newcommand{\DD}{\opname{DD}}
\newcommand{\Drun}{\mathbf{D}}
\renewcommand{\d}{\mathrm{d}}
\newcommand{\Cont}{\mathbf{C}}
\renewcommand{\c}{\mathrm{c}}

\renewcommand{\t}{\mathrm{t}}

\newcommand{\Type}{\mathbf{T}}
\renewcommand{\t}{\mathrm{t}}
\renewcommand{\L}{\opname{d}}

\title{A bijection between Tamari intervals and extended fighting fish}

\author[1]{Enrica Duchi\fnref{fn2}}
\ead{duchi@irif.fr}
\author[1]{Corentin Henriet\corref{cor1}}
\ead{henriet@irif.fr}

\cortext[cor1]{Corresponding author}
\fntext[fn2]{Enrica Duchi was partially supported by grant ANR Combiné}

\address[1]{Institut de Recherche en Informatique Fondamentale, Université Paris Cité, France}

\begin{document}

\begin{abstract}
	
We introduce \emph{extended fighting fish} as branching surfaces that can also be seen as walks in the quarter plane defined by simple rewriting rules. The main result we present in the article is a direct bijection between extended fighting fish and intervals of the Tamari lattice that exchanges multiple natural statistics. The model includes the recently introduced fighting fish of (Duchi, Guerrini, Rinaldi, Schaeffer 2017) that were shown to be equinumerated with synchronized Tamari intervals.

Using the dual surface/walk points of view on extended fighting fish, we show that the area statistic on these fish corresponds to the distance statistic (or maximal length of a chain) in Tamari invervals. We also show that the average area of a uniform random extended fighting fish of size $n$, and hence the average distance over the set of Tamari intervals of size $n$, is of order $n^{5/4}$, in accordance with earlier results for the subclass of ordinary fighting fish.
\end{abstract}

\maketitle

\section{Introduction}\label{sec1}
Fighting fish is a relatively new class of combinatorial objects that has been introduced and enumerated in~\cite{ff} by Duchi, Guerrini, Rinaldi and Schaeffer. Roughly speaking, a fighting fish is a branched surface that is obtained by gluing together flexible unit squares along their edges in a directed way, resulting in independent branches that  may overlap. Fighting fish are counted by  $\frac{2}{(n+1)(2n+1)}\binom{3n}{n}$ like other classical objects as \emph{rooted non-separable planar maps}~\cite{Brown,FPR}, \emph{two-stack sortable permutations}~\cite{Bona,goulden,West,Zeilb}, and \emph{left ternary trees}~\cite{DDP,JS}. The article \cite{ff2} gave a new decomposition of fighting fish extending the classical wasp-waist decomposition of polyominoes and proved that the number of fighting fish with $i$ left lower free edges and $j$ right
lower free edges is $\frac{1}{ij}{2i+j-2\choose j-1}{2j+i-2\choose i-1}$, confirming the apparently close relation of fighting fish with the already cited combinatorial structures. In particular, the authors of~\cite{ff2} proved analytically that fighting fish and left-ternary trees share the same formula with respect to one additional parameter,  the {\em fin length} for the fish and to the {\em core} for the trees. However, to the best of our knowledge, there is currently no known bijection involving these two classes that would explain this equidistribution.

In fact, until now, the only known bijection involving fighting fish was the recursive one given by Fang in~\cite{FangFish} with two-stack sortable permutations, obtained using a new recursive decomposition of this permutation class, which is isomorphic to the one of fighting fish from~\cite{ff2}.
It connects fighting fish bijectively to the other combinatorial structures mentioned above, but only indirectly because this recursive structure is different from the ones used in previously known bijections~\cite{Bona,DDP,FPR,goulden,JS}.

In this paper we present a generalization of fighting fish called
\emph{extended fighting fish}, defined as a branched surface made of
unit squares and triangles, and we give a direct bijection between
these objects and intervals of the Tamari lattices, which specializes to a direct bijection between fighting fish and synchronized
intervals. These results were announced in the extended abstract~\cite{FPSAC22}, together with direct bijections with map families: indeed intervals of the
Tamari lattice of order $n$ have been shown in~\cite{cha06} to be enumerated by
$\frac{1}{(n+1)(2n+1)}\binom{4n+2}{n}$, and these numbers also count \emph{rooted bridgeless planar maps}~\cite{walsh75}, and \emph{rooted simple triangulations}~\cite{tutte62}. These relations to maps, announced in~\cite{FPSAC22}, will be discussed in a subsequent paper and we concentrate here on the link between fighting fish and Tamari lattices.

The bijection presented in this paper gives us an interpretation of the area of extended fighting fish in terms of distance on Tamari intervals, that is, the maximal length of a strictly increasing chain in the interval. The Tamari distance seems to be of importance in the field of trivariate diagonal harmonics~\cite{bergeron2011higher}, and we hope that the formula we obtain in Subsection~\ref{subsec4.2} can help deepen some of the conjectures of~\cite{bergeron2011higher}.

\section{Extended fighting fish}\label{sec2}

\subsection{Fighting fish and extended fighting fish}\label{subsec2.1}

Fighting fish were first defined in 2016 in \cite{ff}. 
as a generalization of parallelogram polyominoes. They are constructed as branching surfaces: unit squares are glued together along their edges. We present in this section extended fighting fish, which is a generalization of fighting fish obtained by attaching to them some additional triangles.

\begin{figure}[H]
	\centering
	\includegraphics[page=1,scale=0.6]{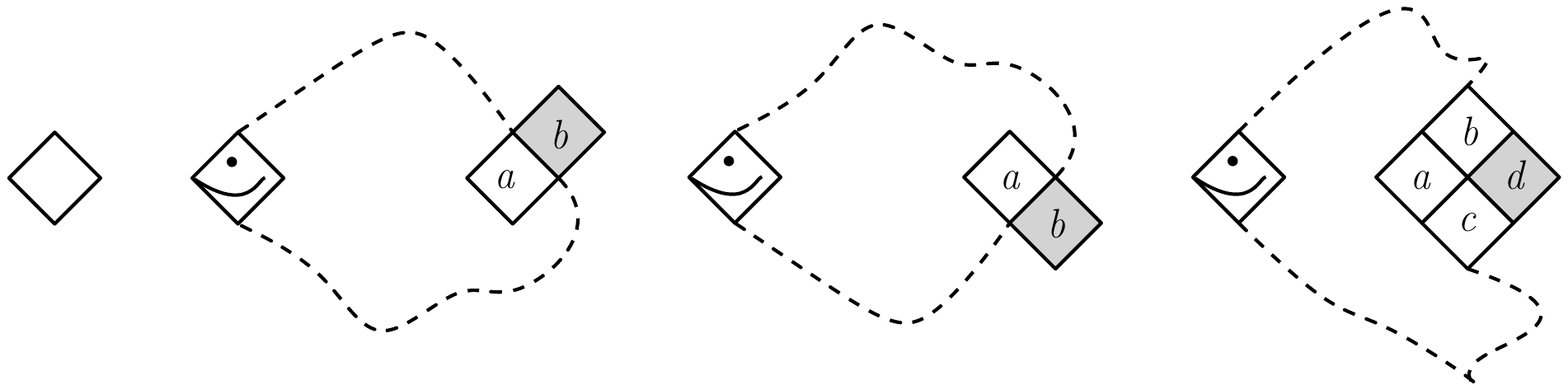}
	\caption{A cell and operations of upper, lower and double gluing.}
	\label{fig:gluings}
\end{figure}

A \emph{cell} is a 45-degree-tilted unit square. Each such cell has a boundary made out of four edges that we call \emph{left lower edge}, \emph{right lower edge}, \emph{right upper edge} and \emph{left upper edge}. We intend to build fighting fish as sets of cells glued together along some of their edges, so we define an edge of a cell to be \emph{free} if it is not glued to the edge of another cell. A \emph{fighting fish} is a finite set of cells constructed by starting with an initial cell (the \emph{head}), then by attaching new cells one by one using one of the three following operations (illustrated in Figure~\ref{fig:gluings}):
\begin{itemize}
	\item \emph{Upper gluing}: Let $a$ be a cell in a fish whose right upper edge is free; we glue the left lower edge of a new cell $b$ to the right upper edge of $a$.
	\item \emph{Lower gluing}: Let $a$ be a cell in a fish whose right lower edge is free; we glue the left upper edge of a new cell $b$ to the right lower edge of $a$.
	\item \emph{Double gluing}: Let $a$, $b$ and $c$ three cells in a fish such that $b$ (resp. $c$) has its left lower (resp. upper) edge glued to the right upper (resp. lower) edge of $a$, and the right lower (resp. upper) edge of $b$ (resp. $c$) is free; we glue both the left upper and left lower edges of a new cell $d$ respectively to the right lower edge of $b$ and to the right upper edge of $c$.
\end{itemize}

While the description of fighting fish is iterative, we are interested in these objects independently of the order in which they are constructed. There can then be multiple ways to grow a given fighting fish. We also want to emphasize that fighting fish are not planar objects in the sense that we cannot always fit them in the plane because some unit squares would represent two or more different cells. Still, we will present them in our two-dimensional pictures by taking care of showing which cells are glued together (Figure~\ref{fig:ffexample} provides an example of two possible representations of the same fighting fish).

\begin{figure}[H]
	\centering
	\includegraphics[page=17,scale=0.7]{figures-these}
	\caption{Two representations of the same fighting fish of size 5 and area 5.}
	\label{fig:ffexample}
\end{figure}

Another way to think about fighting fish is to perform its counterclockwise tour: if we follow the boundary of a fighting fish counterclockwise, starting from the leftmost point of its initial cell (the \emph{nose}), we encounter all its free edges once upon getting back to the nose. We can then state an alternative definition of fighting fish in terms of words. Let us encode each type of free edge by a letter in $\{E, N, W, S\}$: $E$ (resp. $S$) for a left lower (resp. upper) one and $N$ (resp. $W$) for a right lower (resp. upper) one. The set of fighting fish is then the set of finite words on the alphabet $\{E, N, W, S\}$ that can be obtained from the word $ENWS$ using the three operations:
\begin{itemize}
	\item Upper gluing: replace a subword $W$ by $NWS$.
	\item Lower gluing: replace a subword $N$ by $ENW$.
	\item Double gluing: replace a subword $WN$ by $NW$.
\end{itemize}

We can also see these words as excursions (walks starting and ending at the origin) on the square lattice confined to the quarter plane $\{x,y\geq 0\}$ by considering each letter of $\{E, N, W, S\}$ as a unit step in the cardinal direction it denotes.

We define an \emph{unfilled branch point} of a fighting fish to be an occurrence of the subword $WN$ in it. We can now define the class of extended fighting fish in the word setting:

\begin{definition}
	An \emph{extended fighting fish} is a word on the alphabet $\{E, N, W, S, V\}$ obtained from a fighting fish with distinguished unfilled branch points by replacing each marked $WN$ subword by $V$ (see Figure~\ref{fig:fhexample}).
\end{definition}

\begin{figure}[H]
	\centering
	\includegraphics[page=16,scale=0.6]{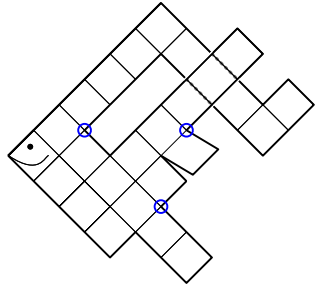}
	\caption{An example of a fighting fish of size 22 and area 26 with its unfilled branch points circled.}
	\label{fig:branchpoints}
\end{figure}

Every $WN \leftrightarrow V$ replacement is clearly revertible so an extended fighting fish can be obtained from a unique fighting fish. We can alternatively see extended fighting fish as walks on $\mathbb{N}^2$ starting and ending at the origin if we consider $V$ as a step $(-1,1)$. In the cell setting, this operation can be seen as the gluing of a triangle instead of a new cell in the operation of double gluing. Fighting fish are then extended fighting fish with no triangle. The edge on the right of a triangle is not declared to be free.

\begin{figure}[H]
	\centering
	\includegraphics[page=18,scale=0.6]{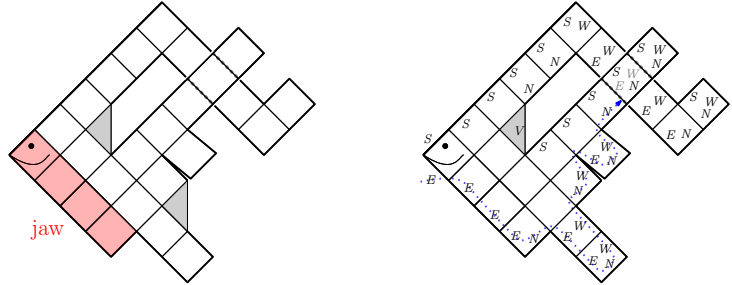}
	\caption{Two of the eight extended fighting fish that can be obtained from the fighting fish of Figure~\ref{fig:branchpoints}.}
	\label{fig:fhexample}
\end{figure}

We define the \emph{size} of an extended fighting fish to be its number of free lower edges minus 1. Every extended fighting fish has as many free lower edges as free upper edges, because this property is true for the head and each type of gluing adds as many free upper edges as free lower edges. Hence we can express the size in the word setting:
$$\size(\FF{F}) = |\opname{F}|_E+|\FF{F}|_N-1 = |\FF{F}|_W+|\FF{F}|_S-1 = \frac{1}{2} \big(|\FF{F}|_E+|\FF{F}|_N+|\FF{F}|_W+|\FF{F}|_S\big) -1\,.$$

We will denote by $\mathcal{EF}_n$ (resp. $\mathcal{F}_n$) the set of extended fighting fish (resp. fighting fish) of size $n$, and $\mathcal{EF} = \bigcup_{n\geq0} \mathcal{EF}_n$ (resp. $\mathcal{F} = \bigcup_{n\geq0} \mathcal{F}_n$). We point out that we need the unions to begin at $n=0$ because we will introduce a special fighting fish of size $0$ in Subsection \ref{subsec2.2}.\\

We also define the \emph{area} of an extended fighting fish to be the number of cells composing it (not counting triangles). It is exactly the algebraic area enclosed by the excursion on the square lattice corresponding to the underlying fighting fish. This observation allows us to give a formula for the area in the word setting analogous to the shoelace formula for polygons. For a word $w \in \{E,N,W,S,V\}^*$, we define its \emph{latitude} and \emph{longitude} to be the $x$ and $y$ coordinates of the endpoint of its corresponding walk on $\mathbb{N}^2$ starting from $(0,0)$. We then have:
\begin{align*}
\lat(w) &= |w|_N - |w|_S + |w|_V\,.\\
\Long(w) &= |w|_E - |w|_W - |w|_V\,.
\end{align*}

Let $\FF{F}$ be an extended fighting fish, $m$ be its length (as a word), and, for $0 \leq j \leq m$, let $\FF{F}^{\leq j}$ be the prefix of $\FF{F}$ of length $j$. Then the area of $\FF{F}$ can be written as:
\[\Area(\FF{F}) = \sum_{i=1}^{m} \Big(\lat\big(\FF{F}^{\leq i}\big)-\lat\big(\FF{F}^{\leq i-1}\big)\Big) \Long\big(\FF{F}^{\leq i}\big)\,.\]

This formula is true for every extended fighting fish because it is true for the head alone and that the left and right-hand sides behave in the same way under the operations of gluing.

\begin{figure}[H]
	\centering
	\includegraphics[page=3,scale=0.6]{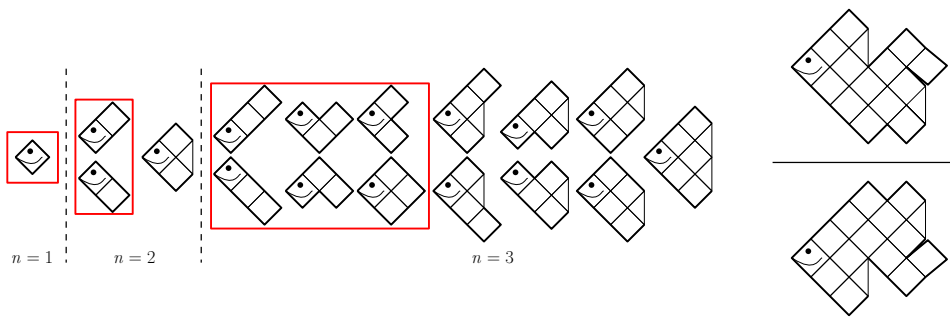}
	\caption{All extended fighting fish of size $n$ for $n = 1,2,3$, with fighting fish in red boxes, and the action of conjugation.}
	\label{fig:smallfish}
\end{figure}

We end this presentation by pointing out that the definitions of extended fighting fish are symmetric with respect to the horizontal axis. Precisely, if we take an extended fighting fish $\FF{F}$ described as a word in $\{E, N, W, S, V\}^*$, reverse it and change its letters with the rules $E \leftrightarrow S$, $N \leftrightarrow W$ (leaving $V$ unchanged), we obtain an extended fighting fish that is the image of $\FF{F}$ under the symmetry with respect to the horizontal axis. We call this extended fighting fish the \emph{conjugate} of $\FF{F}$ and we denote it by $\ConjF(\FF{F})$ (see Figure~\ref{fig:smallfish}). The map $\ConjF$ is then an involution of $\mathcal{EF}$ and of $\mathcal{F}$ that preserves both size and area.

\subsection{Recursive decomposition}\label{subsec2.2}

The wasp-waist decomposition of fighting fish has been presented in~\cite{ff2} and used to present a bijection with two-stack sortable permutations in~\cite{FangFish}. This decomposition relies on the removal of a certain strip at the bottom of the fish and the new decomposition we will present here is in the same spirit, the difference being on the nature of the strip removed. Even if this decomposition can be stated precisely with the word description of extended fighting fish, we prefer here to present its version on gluings of cells that we find more enlightening.\\

When we perform the counterclockwise tour of an extended fighting fish, we first encounter a sequence of free left lower edges before moving on to a free right lower edge. The \emph{jaw} of an extended fighting fish is the strip of cells having their free lower left edge involved in this initial sequence. The \emph{jaw length} is then the number of cells in the jaw, denoted by $\ljaw(\FF{F})$ for an extended fighting fish $\FF{F}$. A \emph{pointed extended fighting fish} is an extended fighting fish $\FF{F}$ having its first $i$ free left lower edges distinguished, for some $1 \leq i \leq \ljaw(\FF{F})+1$. We denote it by $\FF{F}^{\bullet i}$. If $F$ is a fighting fish, $\FF{F}^{\bullet i}$ is said to be \emph{properly pointed} if $i \leq \ljaw(\FF{F})$, that is, we forbid the pointing of the last and only right lower free edge. We denote by $\mathcal{EF}^\bullet_n$ (resp. $\mathcal{F}^\bullet_n$) the set of pointed extended fighting fish (resp. properly pointed fighting fish) of size $n$.\\

Let us now introduce two operations on extended fighting fish:
\begin{itemize}
	
	\begin{figure}[H]
	\centering
	\includegraphics[page=6,scale=0.6]{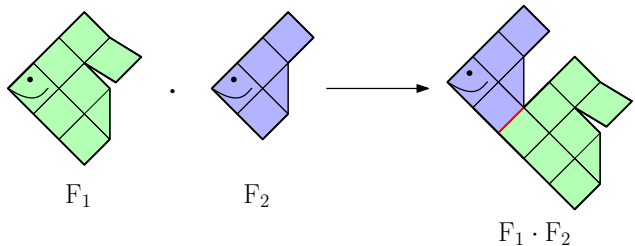}
	\caption{Concatenation of extended fighting fish.}
	\label{fig:fhconc}
	\end{figure}

	\item The \emph{concatenation} $\FF{F}_1 \odot \FF{F}_2$ (see Figure~\ref{fig:fhconc}) of two extended fighting fish $\FF{F}_1$ and $\FF{F}_2$ is obtained by gluing the left upper free edge of the head of $\FF{F}_1$ to the right lower free edge of the last cell of the jaw of $\FF{F}_2$. In the word setting, if we write $\FF{F}_1 = \FF{G}_1 S$ and $\FF{F}_2 = E^{\ljaw(\FF{F}_2)} N \FF{G}_2$, then $\FF{F}_1 \odot \FF{F}_2 = E^{\ljaw(\FF{F}_2)} \FF{G}_1 \FF{G}_2$.
	
\begin{figure}[H]
	\centering
	\includegraphics[page=7,scale=0.75]{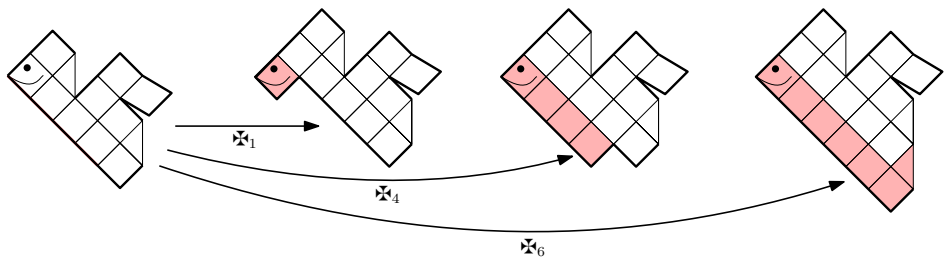}
	\caption{Some augmentations of an extended fighting fish.}
	\label{fig:fhaug}
\end{figure}
	
	\item The \emph{$i$-augmentation} $\maltese_i(\FF{F})$ of an extended fighting fish $\FF{F}$ (written $\FF{F} = E^{\ljaw(\FF{F})} N \FF{G}$ as a word), for $i$ being an integer between 1 and $\ljaw(\FF{F}) + 1$, is obtained in the following way (see Figure~\ref{fig:fhaug}):
	\begin{itemize}
		\item If $i \leq \ljaw(\FF{F})$, then we glue the left lower free edge of each one of the first $i$ cells of the jaw of $\FF{F}$ to a new cell and then glue the right lower edge and the left upper edge of all pairs of adjacent new cells. In the word setting, $\maltese_i(\FF{F}) = E^i N E^{\ljaw(\FF{F})-i} N \FF{G} S$.
		\item If $i = \ljaw(\FF{F}) +1$, then we perform the $\ljaw(\FF{F})$-augmentation on $F$, and we glue the only remaining free right lower edge of the added cells to the left upper edge of a new cell and glue a triangle to the right upper edge of this new cell and to the right lower edge of the final cell of the jaw of $\FF{F}$. In the word setting, $\maltese_{\ljaw(\FF{F})+1}(\FF{F}) = E^{\ljaw(\FF{F})+1} V \FF{G} S$ (while not very natural, adding $V$ in the $(\ljaw(\FF{F}) +1)$-augmentation of $\FF{F}$ is necessary for our following decomposition).
	\end{itemize}
\end{itemize}

These two operations produce valid extended fighting fish. It is indeed not hard to see for the concatenation of $\FF{F}_1$ and $\FF{F}_2$: to construct $\FF{F}_1 \odot \FF{F}_2$, we first build $\FF{F}_2$, then attach a new cell by its left upper edge to the right lower edge of the final cell of the jaw of $\FF{F}_2$ and we finally grow $\FF{F}_1$ starting from this added cell. For the $i$-augmentation of $\FF{F}$, we first remark that every cell of the jaw of $\FF{F}$, except the head, was obtained using a lower gluing because its left lower edge is free. We then construct $\maltese_i(\FF{F})$ by starting from a new head to which we attach the head of $\FF{F}$, and we grow the extended fighting fish like $\FF{F}$ from it upon inserting each new cell just before we insert the jaw cell it will be glued to (so that each one of the $i-1$ first lower gluings of the jaw becomes a right gluing followed by a double gluing). In the case $i = \ljaw(\FF{F})+1$, we additionnally perform a lower gluing on the last cell of the jaw and a triangle gluing. These two operations on extended fighting fish allow us to follow the statistics:

\begin{prop}\label{statsfish}
	For every extended fighting fish $\FF{F}_1$ and $\FF{F}_2$, and every $1 \leq i \leq \ljaw(\FF{F}_1)+1$, the following equalities hold:
	\begin{align*}
	&\size(\FF{F}_1 \odot \FF{F}_2) = \size(\FF{F}_1) + \size(\FF{F}_2)\,, &&\size(\maltese_i(\FF{F}_1)) = \size(\FF{F}_1)+1\,,\\
	&\Area(\FF{F}_1 \odot \FF{F}_2) = \Area(\FF{F}_1) + \Area(\FF{F}_2)\,, &&\Area(\maltese_i(\FF{F}_1)) = \Area(\FF{F}_1) + i\,,\\
	&\ljaw(\FF{F}_1 \odot \FF{F}_2) = \ljaw(\FF{F}_1) + \ljaw(\FF{F}_2)\,, &&\ljaw(\maltese_i(\FF{F}_1)) = i\,.
	\end{align*}
\end{prop}

In order to state our decomposition in an homogeneous way, we define what we call the \emph{empty fish}: it is the word $EVS$, and we can see it as a single triangle. We will refer to it as $\varepsilon$, it has size, area and jaw length all equal to $0$, and we set it to be a fighting fish. We define its concatenations $\FF{F} \odot\hspace{1.8pt}\varepsilon$ and $\varepsilon \odot \FF{F}$ with an extended fighting fish $\FF{F}$ to be both equal to $\FF{F}$ (so $\varepsilon \odot \varepsilon = \varepsilon$) and its $1$-augmentation to be equal to $ENWS$.

\begin{figure}[H]
	\centering
	\includegraphics[page=4,scale=0.75]{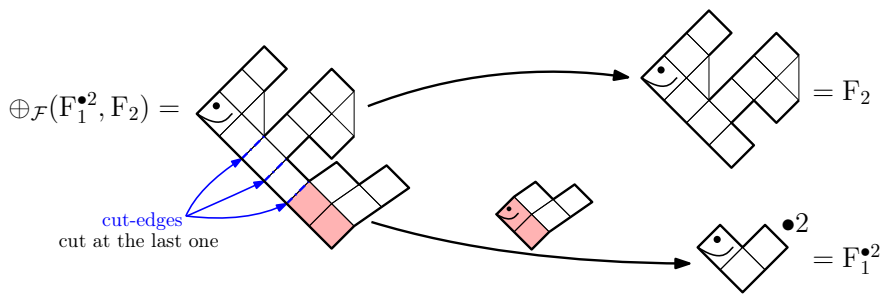}
	\caption{Decomposition of an extended fighting fish.}
	\label{fig:fhdec}
\end{figure}

\begin{thm}\label{decfish}
	Let $\FF{F}_1^{\bullet i}$ be a pointed extended fighting fish of size $n_1$ and $F_2$ be an extended fighting fish of size $n_2$. We define the \emph{low composition} of $\FF{F}_1^{\bullet i}$ and $\FF{F}_2$ to be the extended fighting fish of size $n_1+n_2+1$ equal to:
	$$\oplusF(\FF{F}_1^{\bullet i},\FF{F}_2) = \maltese_i(\FF{F}_1) \odot \FF{F}_2\,.$$
	Then every non-empty extended fighting fish $\FF{F}$ can be decomposed in a unique way as $\FF{F} = \oplusF(\FF{F}_1^{\bullet i},\FF{F}_2)$ (see Figure~\ref{fig:fhdec} for an example), with $\FF{F}_1^{\bullet i}$ being a pointed extended fighting fish and $\FF{F}_2$ being an extended fighting fish, and $\FF{F}$ has jaw length $\ljaw(\FF{F}_2) + i$ and area $\Area(\FF{F}_1) + \Area(\FF{F}_2) + i$. In other words, $\oplusF$ is a bijection from $\mathcal{EF}^\bullet \times \mathcal{EF}$ to $\mathcal{EF}-\{\varepsilon\}$. Furthermore, $\oplusF$ induces a bijection from $\mathcal{F}^\bullet \times \mathcal{F}$ to $\mathcal{F}-\{\varepsilon\}$.
\end{thm}

\begin{proof}
	We already know that low composition produces valid extended fighting fish. Let us now prove that the decomposition exists and is unique. Let $\FF{F}$ be a non-empty extended fighting fish having jaw length $k \geq 1$. If $\FF{F}$ is the single cell extended fighting fish, we know that it has the unique decomposition $\oplusF(\varepsilon^{\bullet 1}, \varepsilon)$, so we now take $\FF{F}$ as different from $ENWS$. We define a \emph{cut-edge} of $\FF{F}$ to be a common edge between two adjacent cells of its jaw such that ungluing this edge separates $\FF{F}$ in two connected components (see Figure~\ref{fig:fhdec}). Let $a$ and $b$ be two adjacent cells in the jaw such that the right lower edge of $a$ is glued to the left upper edge of $b$. Then the common edge between $a$ and $b$ is never a cut-edge if a triangle is glued to the right upper edge of $b$ (because in this case the left upper edge of this triangle is glued to the right lower edge of another cell $c$ which is glued by its left lower edge to the right upper edge of $a$), and always one if a triangle is glued on the right upper edge of $a$. Otherwise it is a cut-edge if and only if the cells glued to the right upper edge of $a$ and $b$ are not glued together. Two situations are to be distinguished:
	\begin{itemize}
		\item If the jaw of $\FF{F}$ contains no cut-edge, then removing all cells of its jaw (along with the possible triangle if one is glued on the right upper edge of the last cell of the jaw) yields a non-empty extended fighting fish $\FF{F}_1$ having jaw length at least $k-1$. Then $\FF{F}$ is just the $k$-augmentation of $\FF{F}_1$, that is, $\FF{F} = \oplusF(\FF{F}_1^{\bullet k}, \varepsilon)$.
		\item If the jaw of $\FF{F}$ contains at least one cut-edge, then cutting $\FF{F}$ at one of them yields two connected components $\tilde{\FF{F}}_1$ and $\FF{F}_2$ that are non-empty extended fighting fish and such that $\FF{F} = \tilde{\FF{F}}_1 \odot \FF{F}_2$. The only cut-edge that can be cut in order to have $\tilde{\FF{F}}_1$ without cut-edge is the last one in the jaw (the furthest from the head). Then, when cutting this way, thanks to the previous case, $\tilde{\FF{F}}_1$ can be written as $\maltese_i(\FF{F}_1)$ with $\FF{F}_1$ being a (possibly empty) extended fighting fish, and then $\FF{F} = \oplusF(\FF{F}_1^{\bullet i}, \FF{F}_2)$.
	\end{itemize}
	
	The decomposition is unique because the $\maltese$ operation is clearly revertible, the concatenation of two non-empty extended fighting fish yields at least one cut-edge and in this case, the cut-edge where we can cut is unique. Correspondence of statistics are derived from Proposition~\ref{statsfish} and the restriction to fighting fish comes easily by the observation that adding a triangle is equivalent to applying the operation $\maltese$ with the maximal index.
\end{proof}

Let us point out that we chose to cut the last cut-edge to decompose an extended fighting fish with respect to $\oplusF$. In fact, if we choose to decompose an extended fighting fish into two smaller extended fighting fish by cutting it at the first cut-edge, it gives rise to another decomposition, isomorphic to the first one we presented. The proof is similar to the one above.

\begin{prop}
		Let $\FF{F}_1^{\bullet i}$ be a pointed extended fighting fish of size $n_1$ and $\FF{F}_2$ be an extended fighting fish of size $n_2$. We define the \emph{high composition} of $\FF{F}_1^{\bullet i}$ and $\FF{F}_2$ to be the extended fighting fish of size $n_1+n_2+1$, jaw length $\ljaw(\FF{F}_2)+i$ and area $\Area(\FF{F}_1) + \Area(\FF{F}_2) + i$ defined by:
		$$\otimesF(\FF{F}_1^{\bullet i},\FF{F}_2) = \FF{F}_2 \odot \maltese_i(\FF{F}_1)\,.$$
		Then $\otimesF$ is a bijection from $\mathcal{EF}^\bullet \times \mathcal{EF}$ to $\mathcal{EF}-\{\varepsilon\}$. Furthermore, $\otimesF$ induces a bijection from $\mathcal{F}^\bullet \times \mathcal{F}$ to $\mathcal{F}-\{\varepsilon\}$.
\end{prop}

\section{Tamari intervals}\label{sec3}

The Tamari lattice have been introduced by Dov Tamari in his seminal work presented in \cite{tamarithesis}. It gives a structure of order to the set of Dyck paths of a given size. The enumeration of the intervals of the Tamari lattice has first been performed by Chapoton in~\cite{cha06}, and following this work, the subclass of synchronized intervals was further enumerated by Fang and Préville-Ratelle in~\cite{FPR}. The purpose of this section is to present a decomposition of Tamari intervals that was stated by Bousquet-Mélou, Fusy and Préville-Ratelle in~\cite{BMFPR11}, which is isomorphic to the one described for extended fighting fish in the previous section. We will also follow some statistics on the decomposition thanks to results obtained by Pons in~\cite{Pons19} in the setting of interval-posets.


\subsection{The Tamari lattice on Dyck paths}\label{subsec3.1}

\begin{definition}
	A \emph{Dyck path} of size $n$ is a finite walk on $\mathbb{Z}^2$ starting at (0,0) that consists of $n$ up steps $u=(1,1)$, $n$ down steps $d=(1,-1)$ and that stays above the $x$-axis.	A \emph{peak} of $P$ is a point of $P$ preceded by an up step and followed by a down step. A \emph{valley} of $P$ is a point of $P$ preceded by a down step and followed by an up step.	A \emph{double rise} (resp. \emph{double descent}) of $P$ is a point of $P$ between two up steps (resp. between two down steps).	A \emph{contact} of $P$ is a point of $P$ lying on the $x$-axis.
\end{definition}

Note that every non-empty Dyck path $P$ can be decomposed uniquely $P = P_1 u P_2 d$ with $P_1$ and $P_2$ being Dyck paths, by cutting it in two parts at its first return to the $x$-axis.\\
We will denote by $\mathcal{D}_n$ the set of Dyck paths of size $n$, by $\mathcal{D} = \bigcup_{n \geq 0} \mathcal{D}_n$ the set of all Dyck paths, and, for P a Dyck path, by $\Val(P)$, $\Peak(P)$, $\DR(P)$ and $\DD(P)$ the number of valleys (occurences of $du$), of peaks ($ud$), of double rises ($uu$) and of double descents ($dd$) of $P$ respectively. Note that we always have $\Val(P)+1 = \Peak(P)$, $\DR(P)=\DD(P)$ and $\Val(P)+\DR(P)=n-1$.\\

\begin{figure}[H]
	\centering
	\includegraphics[page=10,scale=0.6]{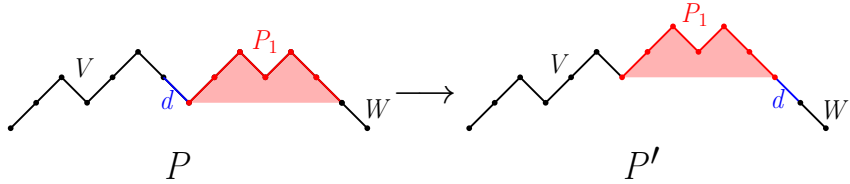}
	\caption{Example of right rotation on a Dyck path.}
	\label{fig:rightrotation}
\end{figure}

We define a covering relation on $\mathcal{D}_n$. Let $P$ be an element of $\mathcal{D}_n$ that can be written (as a word) as $P = V d P_1 W$, where $V,W$ are words on the alphabet $\{u,d\}$ and $P_1$ is a non-empty Dyck path that returns to the $x$-axis only at the end. We then construct $P' = V P_1 d W$, which is also a Dyck path, and we say that $P'$ covers $P$ and that $P'$ can be obtained from $P$ by \emph{right rotation} (see Figure~\ref{fig:rightrotation}). The \emph{Tamari lattice} $(\mathcal{D}_n, \preceq)$ of order $n$ is given by the transitive closure $\preceq$ of the covering relation we just defined: $P \preceq P'$ if $P'$ can be obtained from $P$ through a (possibly empty) series of right rotations. 

We recall now the conjugation of Dyck paths. We refer to Subsection $4.4$ of~\cite{FMN} for an extended version of results presented here. Conjugation of Dyck paths is an anti-automorphism $\ConjD$ of Tamari lattices, that is, a bijection from $\D_n$ to $\D_n$ for every $n \geq 0$ such that for all $P,Q \in \D_n$, we have $P \preceq Q$ if and only if $\ConjD(Q) \preceq \ConjD(P)$. Let us define recursively $\ConjD(\bullet) = \bullet$ and for all Dyck paths $P_1$ and $P_2$, $\ConjD(P_1 u P_2 d) = \ConjD(P_2) u \ConjD(P_1) d$.

\begin{prop}[\cite{FMN}, Proposition 4.7]
	$\ConjD$ is an anti-automorphism and an involution of the Tamari lattice $\D_n$.
\end{prop}

\begin{proof}
	We have $\ConjD^2(\bullet) = \bullet$, and, for all Dyck paths $P_1$ and $P_2$, 
	$$\ConjD^2(P_1uP_2d) = \ConjD(\ConjD(P_2) u \ConjD(P_1) d) = \ConjD^2(P_1) u \ConjD^2(P_2) d\,.$$
	
	It then follows by induction that $\ConjD^2(P) = P$ for every Dyck path, hence $\ConjD$ is an involution. The anti-automorphism property of $\ConjD$ can be proven using the unique decomposition $P =P_1 u P_2 d$ of Dyck paths, but a more enlightening way to visualize it is to see the Tamari lattices as posets on complete binary trees ordered by right rotations: the conjugation corresponds in this setting to recursively switching left and right subtrees. We refer to Subsection 1.1 of~\cite{PRV} for a detailed presentation of the Tamari lattice in terms of complete binary trees.
\end{proof}

Tamari lattices are then self-dual via the operation of conjugation, and we introduce now some vectors that behave nicely with respect to this conjugation.

\begin{figure}[H]
	\centering
	\includegraphics[page=9,scale=0.6]{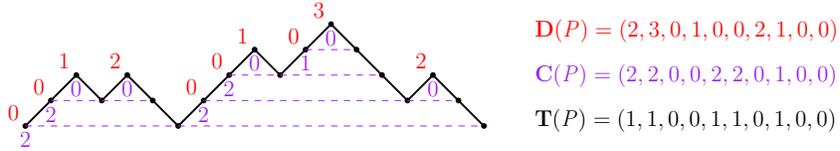}
	\caption{A Dyck path, its descent and contact vectors, and its type.}
	\label{fig:dyck-vectors}
\end{figure}

\begin{definition}(see Figure~\ref{fig:dyck-vectors})
	Let $P$ be a Dyck path of size $n$.
	\begin{itemize}
	\item The \emph{descent vector} of $P$ is the vector of nonnegative integers $\Drun(P)=(\d_0(P),...,\d_n(P))$ such that $P = d^{\d_n(P)}ud^{\d_{n-1}(P)}u...ud^{\d_0(P)}$.
	\item The \emph{contact vector} of $P$ is the vector of nonnegative integers $\Cont(P)=(c_0(P),...,c_n(P))$ such that $c_0(P)$ is the number of non-initial contacts of $P$ and $c_i(P)$ is the number of non-initial contacts of the Dyck path following the $i^{\rm{th}}$ up step of $P$.
	\item The \emph{type} of $P$ is the binary vector $\Type(P) = (\t_0(P),...,\t_n(P))$ where $t_i(P)=\mathds{1}_{c_i(P)>0}$.
	\end{itemize}
\end{definition}

Note that a Dyck path is fully determined by the data of its descent vector (or by its contact vector). Peaks of a Dyck path $P$ correspond exactly to zero components of the contact vector, hence we have:
\begin{prop}\label{valdr}
	\[\Val(P)+1 = |\{0 \leq i \leq n, \c_i(P)=0\}|\,,\]
	\[\DR(P)+1 = |\{0 \leq i \leq n, \c_i(P)\geq 1\}|\,.\]
\end{prop}

Also, conjugation exchanges the descent and the contact vectors:

\begin{prop}
	For every Dyck path $P$, we have $\Cont(P) = \Drun\big(\Conj_\D(P)\big)$.
\end{prop}

\begin{proof}
	We proceed by induction on the size of Dyck paths:
	\begin{itemize}
		\item The empty Dyck path is its own conjugate and has contact and descent vectors both equal to 0.
		\item Let $P$ be a Dyck path of size $n\geq 0$ that we write as its unique decomposition $P=P_1uP_2d$, with $P_1$ and $P_2$ being Dyck paths of respective sizes $n_1$ and $n_2$, strictly smaller than $n$. Using induction hypothesis on $P_1$ and $P_2$, we then have
		\begin{align*}
		\Drun\big(\ConjD(P)\big) &= \Drun\big(\ConjD(P_2)u\ConjD(P_1)d\big)\\
		&= \big(\c_0(P_1)+1, \c_1(P_1), ..., \c_{n_1}(P_1), \c_0(P_2), ..., \c_{n_2}(P_2)\big)\\
		&= \Cont(P)\,.
		\end{align*}
	\end{itemize}
\end{proof}

\subsection{Recursive decomposition of Tamari intervals}\label{subsec3.2}

We now want to recall two isomorphic decompositions of Tamari intervals. These decompositions are also isomorphic to the ones presented for extended fighting fish in Subsection \ref{subsec2.2}, and we will use them to construct a bijection between extended fighting fish and Tamari intervals that specializes nicely to the subclass of fighting fish. 

\begin{figure}[H]
	\centering
	\includegraphics[page=12,scale=0.6]{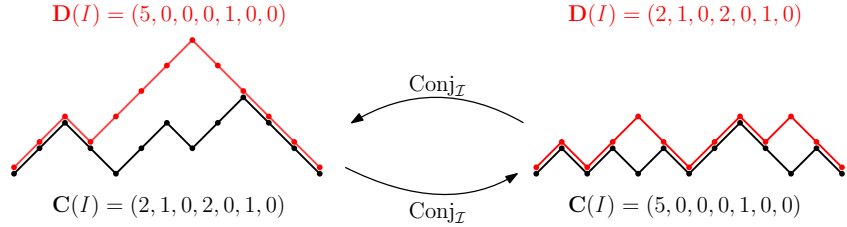}
	\caption{Contact and descent vectors of an interval and the action of conjugation.}
	\label{fig:intconj}
\end{figure}

\begin{definition} (see Figure~\ref{fig:intconj})
	A \emph{Tamari interval} of size $n\geq 0$ is an interval in the Tamari lattice $(\mathcal{D}_n, \preceq)$ of order $n$, that is, a pair of comparable Dyck paths $[P, Q]$ with $P \preceq Q$. For such an interval $I=[P,Q]$, we define:
	\begin{itemize}
		\item its contact vector $\Cont(I)=\Cont(P)$,
		\item its descent vector $\Drun(I)=\Drun(Q)$,
		\item its conjugate $\ConjI(I)=[\ConjD(Q),\ConjD(P)]$, and
		\item its \emph{Tamari distance} $\L(I)$ to be the length of the longest strictly increasing chain from $P$ to $Q$ in the Tamari lattice. A strictly increasing chain of length k from $P$ to $Q$ is a list of Dyck paths $P= D_1 < D_2 < ... < D_k = Q$.
	\end{itemize}

	We say that a Tamari interval is \emph{synchronized} if $P$ and $Q$ have the same type, and we set $\type(I)=\type(P)=\type(Q)$.\\
	
	A \emph{pointed Tamari interval} is an interval $[P,Q]$ with a distinguished contact of the lower path $P$, which we write as $[P^\ell {\cdot} P^r,Q]$, where $P=P^\ell P^r$ is split by the distinguished contact into two subpaths that are also Dyck paths $P^\ell$ and $P^r$ (that can be empty). We also write $[P, Q]^{\bullet i}$ for some $1 \leq i \leq \c_0(I)+1$ to denote the pointed interval obtained from $[P,Q]$ by distinguishing the $i^{\rm{th}}$ contact from right to left of $P$ (while this way of pointing from right to left may be counter-intuitive, it is the right notion to consider for the decomposition to work). \\
	
	A \emph{properly pointed synchronized interval} is a pointed synchronized interval $[P^\ell {\cdot} P^r,Q]$ such that $P^\ell$ is non-empty, or equivalently, it is a pointed synchronized interval $[P, Q]^{\bullet i}$, with $1 \leq i \leq \c_0(I)$.
\end{definition}

We borrow the expression "Tamari distance" from \cite{bergeron2011higher}, but let us point out that it is not a distance in the topological sense, as $\L([P,Q])$ is not the distance between $P$ and $Q$ in the Hasse diagram of the Tamari lattice. However, the quantity $\L([P,Q])$ still expresses how far $Q$ is from $P$ in the Tamari lattice, so we decided to keep the notation of \cite{bergeron2011higher}.\\

We denote by $\mathcal{I}_n$ (resp. $\mathcal{SI}_n$) the set of Tamari intervals (resp. synchronized intervals) of size $n$ and by $\mathcal{I}^\bullet_n$ (resp. $\mathcal{SI}^\bullet_n$) the set of pointed Tamari intervals (resp. properly pointed synchronized intervals) of size $n$.

\begin{figure}[H]
	\centering
	\includegraphics[page=8,scale=0.8]{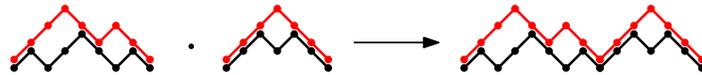}
	\caption{Concatenation of two intervals.}
	\label{fig:intconc}
\end{figure}

\begin{figure}[H]
	\centering
	\includegraphics[page=15,scale=0.75]{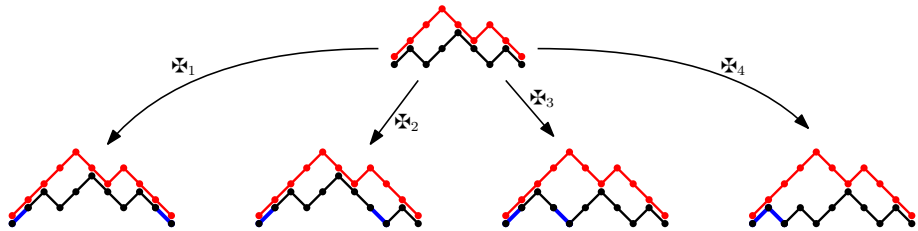}
	\caption{Augmentations of an interval.}
	\label{fig:intaug}
\end{figure}

In order to express decompositions of Tamari intervals, we define on Tamari intervals the following operations of concatenation and augmentation:

\begin{definition}
	Let $I_1 = [P_1, Q_1]$ and $I_2=[P_2,Q_2]$ be two Tamari intervals, and let $1 \leq i \leq \c_0(I_1) +1$ such that $I_1^{\bullet i}=[P_1^\ell \cdot P_1^r,Q_1]$ is a well-defined pointed interval.
	\begin{itemize}
	\item The \emph{concatenation} of $I_1$ and $I_2$ is the Tamari interval defined as $I_1 \odot I_2 = [P_1 P_2, Q_1 Q_2]$ (see Figure~\ref{fig:intconc}).
	\item The \emph{$i$-augmentation} of $I_1$ is the Tamari interval $\maltese_i(I_1) = [u P_1^\ell d P_1^r, u Q_1 d]$ (see Figure~\ref{fig:intaug}).
	\end{itemize}
\end{definition}

These two operations appear in the setting of interval-posets developped in~\cite{Pons19}. Precisely, when translating the definitions appearing in Definition 31 of~\cite{Pons19} in terms of Tamari intervals, our operation of concatenation of $I_1$ and $I_2$ is the same as the operation of left grafting of $I_1$ over $I_2$ and our operation of $i$-augmentation of $I_1$ corresponds to the $r$-right grafting of $I_1$ over $[ud,ud]$, with $r = \c_0(I_1)+1-i$. In Section $4.1$ of~\cite{Pons19}, the behavior of statistics such as Tamari distance with respect to the operations $\odot$ and $\maltese$ is studied, we summarize it in the following proposition:

\begin{prop}[Adapted from~\cite{Pons19}]\label{opstats}
	Let $I_1$ and $I_2$ be two Tamari intervals of respective size $n_1$ and $n_2$, and let $1 \leq i \leq \c_0(I_1)+1$. We have the following relations:
	\begin{align*}
	\Cont(I_1 \odot I_2) &= (\c_0(I_1)+\c_0(I_2),\c_1(I_1),...,\c_{n_1}(I_1),\c_1(I_2),...,\c_{n_2}(I_2))\,,\\
	\Drun(I_1 \odot I_2) &= (\d_0(I_2),...,\d_{n_2-1}(I_2),\d_0(I_1),...,\d_{n_1-1}(I_1),0 = \d_{n_2}(I_2)+\d_{n_1}(I_1))\,,\\
	\L(I_1 \odot I_2) &= \L(I_1) + \L(I_2)\,,
	\end{align*}
	\begin{align*}
	\Cont(\maltese_i(I_1)) &= (i, \c_0(I_1)+1-i,\c_1(I_1),...,\c_{n_1}(I_1))\,,\\
	\Drun(\maltese_i(I_1)) &= (\d_0(I_1)+1,\d_1(I_1),...,\d_{n_1}(I_1),0)\,,\\
	\L(\maltese_i(I_1)) &= \L(I_1) +i-1\,.
	\end{align*}
\end{prop}

%

We point out that although the behavior of contact and descent vectors is easy to derive from the definitions of operations, the relations concerning the Tamari distance are much less evident. The proof of~\cite{Pons19} is quite involved and relies on the fact that performing allowed right rotations from left to right to get to $Q$ from $P$ give rise to a longest strictly increasing chain from $P$ to $Q$.

We are now able to state the decomposition of Tamari intervals that appears in~\cite{BMFPR11} in the general setting of $m$-Tamari lattices (here $m=1$). This decomposition specializes nicely to a decomposition of synchronized intervals, that is exactly the one that appears in~\cite{FPR}.

\begin{figure}[H]
	\centering
	\includegraphics[page=13,scale=0.9]{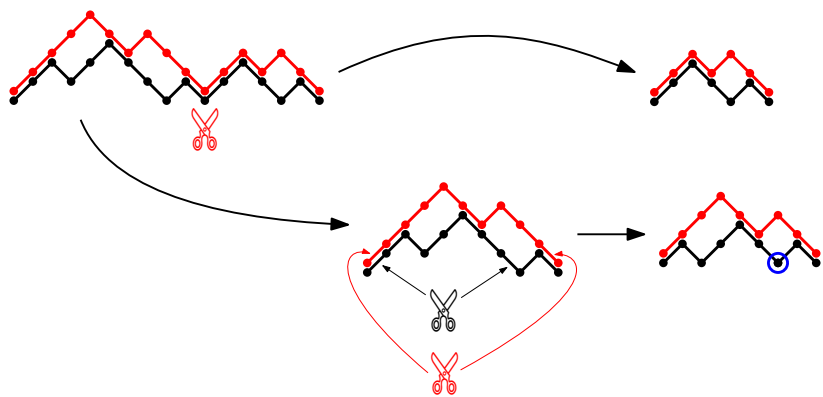}
	\caption{Decomposition of Tamari intervals.}
	\label{fig:int-decomp}
\end{figure}

\begin{prop}[\cite{BMFPR11}, Section 2.3]\label{decint}
	Let $I^{\bullet i}_1 = [P_1^\ell {\cdot} P_1^r,Q_1]$ be a pointed Tamari interval of size $n_1$ with $1\leq i \leq \c_0(I_1)+1$ and $I_2 = [P_2, Q_2]$ be a Tamari interval of size $n_2$. The \emph{left composition} of $I^{\bullet i}_1$ and $I_2$ is defined to be the Tamari interval of size $n_1+n_2+1$  (see Figure~\ref{fig:int-decomp}):
	$$\oplusI(I^{\bullet i}_1, I_2) = \maltese_i(I_1) \odot I_2\,.$$
	Then the map $\oplusI$ is a bijection from $\mathcal{I}^\bullet \times \mathcal{I}$ to $\mathcal{I} - \{[\bullet,\bullet]\}$.
\end{prop}

\begin{prop}[\cite{FPR}, Proposition 3.2]
	The map $\oplusI$ induces a bijection between $\mathcal{SI}^\bullet \times \mathcal{SI}$ and $\mathcal{SI} - \{[\bullet,\bullet]\}$.
\end{prop}

A simple composition of the relations of Proposition~\ref{opstats} allows us to compute the statistics of $\oplusI(I_1^{\bullet i},I_2)$:
\begin{prop}\label{statsint}
\begin{align*}
\Cont(\oplusI(I^{\bullet i}_1, I_2)) &= (\c_0(I_2)+i,\c_0(I_1)+1-i,\c_1(I_1),...,\c_{n_1}(I_1),\c_1(I_2),...,\c_{n_2}(I_2))\,,\\
\Drun(\oplusI(I^{\bullet i}_1, I_2)) &= (\d_0(I_2),...,\d_{n_2-1}(I_2),\d_0(I_1)+1,\d_1(I_1),...,\d_{n_1}(I_1),\d_{n_2}(I_2))\,,\\
\L(\oplusI(I^{\bullet i}_1, I_2)) &= \L(I_1) + \L(I_2) + i - 1\,.
\end{align*}
In particular, $\c_0(\oplusI(I^{\bullet i}_1, I_2)) = \c_0(I_2)+i$.\\
\end{prop}

This decomposition of (synchronized) intervals allows us to prove the following fact about the nullity of components of contact and descent vectors:

\begin{prop}\label{compzero}
	Let $I$ be a Tamari interval of size $n$. Then for every $0 \leq i \leq n$, either $\c_i(I)$ or $\d_{n-i}(I)$ is $0$.
\end{prop} 

\begin{proof}
	We prove it by induction on the size of the interval:
	\begin{itemize}
		\item Contact and descent vectors of $[\bullet, \bullet]$ are both equal to $(0)$, so we get the base case.
		\item Let $I$ be a Tamari interval of size $n \geq 1$. By Proposition~\ref{decint}, $I$ can be uniquely decomposed as $I = \oplusI(I^{\bullet i}_1, I_2)$ with $I_1$ and $I_2$ of respective sizes $n_1$ and $n_2$, strictly smaller than $n$. Following the contact and descent vectors thanks to Proposition~\ref{statsint}, all but three of the $\c_i(I)\d_{n-i}(I)$ are obviously 0 because they come either from $I_1$ or from $I_2$. The three last products are $\c_0(I)\d_n(I) = (\c_0(I_2)+i)\d_{n_2}(I_2)$, $\c_1(I)\d_{n-1}(I) = (\c_0(I_1)+1-i)\d_{n_1}(I_1)$ and $\c_{n_1+1}(I)\d_{n_2}(I) = \c_{n_1}(I_1)(\d_0(I_1)+1)$, and they are all equal to zero because the last component of a contact/descent vector is always zero.
	\end{itemize}
\end{proof}

To be complete, let us point out that another decomposition isomorphic to the one presented above appears in Proposition $38$ of~\cite{Pons19}. It consists in cutting the Tamari interval at the last contact of its upper path with the $x$-axis.

\begin{prop}[Adapted from~\cite{Pons19}]
	Let $I_1^{\bullet i}$ be a pointed Tamari interval of size $n_1$ and $I_2$ be a Tamari interval of size $n_2$. We define the \emph{right composition} of $I_1^{\bullet i}$ and $I_2$ to be the Tamari interval of size $n_1+n_2+1$ defined by:
	$$\otimesI(I_1^{\bullet i},I_2) = I_2 \odot \maltese_i(I_1)\,.$$
	Then $\otimesI$ is a bijection from $\mathcal{I}^\bullet \times \mathcal{I}$ to $\mathcal{I}-\{[\bullet, \bullet]\}$. Furthermore, $\otimesI$ induces a bijection from $\mathcal{SI}^\bullet \times \mathcal{SI}$ to $\mathcal{SI}-\{[\bullet, \bullet]\}$.
\end{prop}

Remark that we have the same relations as with $\oplus_\mathcal{I}$ for $\c_0$ and $\L$:
\begin{align*}
\c_0(\otimesI(I^{\bullet i}_1, I_2)) &= \c_0(I_2) + i\,,\\
\L(\otimesI(I^{\bullet i}_1, I_2)) &= \L(I_2) + \L(I_1) + i - 1\,.
\end{align*}

We point out that, because $\oplusI$ and $\otimesI$ are isomorphic decompositions, it is possible to define a distance-preserving involution $\Gamma$ on Tamari intervals by:
\begin{align*}
\Gamma([\bullet,\bullet]) &= [\bullet,\bullet]\,,\\
\Gamma(\otimesI(I^{\bullet i}_1, I_2)) &= \oplusI(\Gamma(I_1)^{\bullet i}, \Gamma(I_2))\,.
\end{align*}

This involution appears in terms of grafting trees in~\cite{Pons19}, where it is called the left-branch involution. In the same article, Pons combine it with the conjugation involution to obtain interesting results of equidistribution of certain statistics on Tamari intervals. 

%
%
%
%
%
%
%

\section{The bijection and applications}\label{sec4}

\subsection{Bijection between Tamari intervals and extended fighting fish}\label{subsec4.1}

In the two preceeding sections, we have seen that the given decompositions of extended fighting fish and Tamari intervals are isomorphic, and this fact enables us to recursively define a bijection $\Phi$ from Tamari intervals to extended fighting fish. We set:
\[\begin{cases*}
\Phi([\bullet, \bullet]) = \varepsilon\,,\\
\Phi(\oplusI(I_1^{\bullet i},I_2)) = \oplusF(\Phi(I_1)^{\bullet i},\Phi(I_2)) \text{ for all } (I_1^{\bullet i},I_2) \in \mathcal{I}^\bullet \times \mathcal{I}\,.
\end{cases*}\]

\begin{thm}\label{recbij}
	$\Phi: \mathcal{I} \rightarrow \mathcal{EFF}$ is a bijection such that for every Tamari interval $I$ of size $n \geq 0$, the extended fighting fish $\Phi(I)$ has size $n$ and the following relations hold:
	\[\ljaw(\Phi(I)) = \c_0(I)\,, \hspace{3cm} \Area(\Phi(I)) = \L(I) + n\,.\]
	Moreover, the restriction of $\Phi$ to synchronized intervals induces a bijection from $\mathcal{SI}$ to $\mathcal{FF}$.
\end{thm}

\begin{proof}
	To see that $\Phi$ is well-defined and bijective, we have to prove by induction on $n \geq 0$ that for every Tamari interval $I$ of size $n$, $\Phi(I)$ has size $n$ and can be pointed by exactly the same integers, that is to say $\c_0(I)=\ljaw(\Phi(I))$. Those conditions are true at every step of the induction because they are true for $I = [\bullet,\bullet]$ and composition relations for size and for $\c_0$ and $\ljaw$ given in Theorems~\ref{decfish} and~\ref{statsint} imply the inductive step. The same reasoning applies to deduce that $\Phi$ induces a bijection from $\mathcal{SI}$ to $\mathcal{EFF}$.\\
	
	The equality $\c_0(I)=\ljaw(\Phi(I))$ is already proved by the induction above. If we denote by $\opname{size}(I)$ the size of an interval $I$, then by Proposition~\ref{statsint}, we have:
	\begin{align*}\L(\oplusI(I^{\bullet i}_1, I_2)) + \opname{size}(\oplusI(I^{\bullet i}_1, I_2)) &= \L(I_1) + \L(I_2) + i - 1 + \opname{size}(I_1) + \opname{size}(I_2) +1\\
	&= \L(I_1) + \opname{size}(I_1) + \L(I_2) + \opname{size}(I_2) + i\,.
	\end{align*}
	
	On the other hand, by Theorem~\ref{decfish}, we have: $\Area(\oplusF(\FF{F}_1^{\bullet i},\FF{F}_2)) = \Area(\FF{F}_1) + \Area(\FF{F}_2) + i$. Hence $\L + \opname{size}$ and $\Area$ satisfy the same composition relations and they are both equal to 0 for the empty structure, so we get the desired equality.
\end{proof}

Using Theorem~\ref{recbij} and the enumeration of Tamari intervals in~\cite{cha06}, we get the following enumeration corollary about extended fighting fish:

\begin{cor}
	The number of extended fighting fish of size $n$ is:
	$$|\mathcal{EF}_n| = \frac{2}{n(n+1)} \binom{4n+1}{n-1}\,.$$
\end{cor}

\begin{figure}[H]
	\centering
	\includegraphics[page=14,scale=0.69]{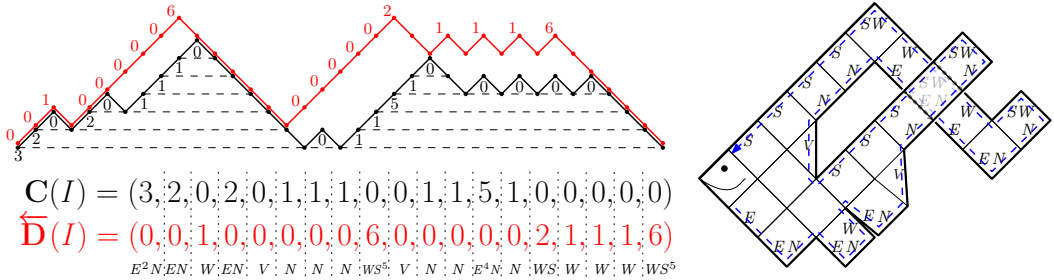}
	\caption{An example of the bijection between Tamari intervals and extended fighting fish.}
	\label{fig:example-bij-tamari}
\end{figure}

The recursive definition of the bijection $\Phi$ is not completely satisfying because it does not provide us a full understanding of the correspondence between the properties of Tamari intervals and extended fighting fish. For example, we will see that the conjugation of Tamari intervals is mapped by $\Phi$ to the symmetry with respect to the horizontal axis for extended fighting fish, but this is not evident with the recursive definition of $\Phi$. Indeed, these notions of symmetries do not behave nicely under the recursive decompositions. We then have to seek for a direct description of $\Phi$ and the recursive description of this bijection will help us find it by enabling us to compute extended fighting fish corresponding to Tamari intervals of small size, allowing us to conjecture the way to de-recursify. We now present such a direct description of $\Phi$ (see Figure~\ref{fig:example-bij-tamari} for an example):

\begin{thm}
	Let $I=[P,Q]$ be a Tamari interval of size $n$, with $\Cont(I)=\Cont(P)$ and $\Drun(I)=\Drun(Q)$ its contact and descent vectors. For $0 \leq i \leq n$, either $\c_i(I)=0$ or $\d_{n-i}(I) = 0$ (see Proposition~\ref{compzero}), and we define $w_i$ by:
	\begin{align*}
	&w_i = E^{\c_i(I)-1}N &\text{ if } \c_i(I) \geq 1 \text{ and } \d_{n-i}(I) = 0\,,\\
	&w_i = WS^{\d_{n-i}(I)-1} &\text{ if } \c_i(I) =0 \text{ and } \d_{n-i}(I) \geq 1\,,\\
	&w_i = V &\text{ if } \c_i(I) =0 \text{ and } \d_{n-i}(I) = 0\,.
	\end{align*}
	Then the word $\FF{F} = Ew_0w_1...w_nS$ is the extended fighting fish $\FF{F} = \Phi(I)$.
\end{thm}

\begin{proof}
	We proceed by induction on the size $n$ of the Tamari interval. For $n = 0$, the interval $[\bullet, \bullet]$ has contact and descent vectors both equal to $(0)$, so $w_0 = V$ and then we get $Ew_0S = EVS = \Phi([\bullet, \bullet])$. We now consider a Tamari interval $I$ of size $n \geq 1$ and assume that we have the equality for every Tamari interval of size strictly less than $n$. We write $I = \oplusI(I_1^{\bullet i}, I_2)$, with $I_1$ of size $n_1$, $I_2$ of size $n_2$ and $n_1+n_2+1 = n$. According to the induction hypothesis, we have $\Phi(I_1) = Ew_0^{(1)}...w_{n_1}^{(1)}S$ and $\Phi(I_2) = Ew_0^{(2)}...w_{n_2}^{(2)}S$, with, for $\alpha \in \{1,2\}$ and $0 \leq j \leq n_\alpha$,
	\begin{align*}
	&w_j^{(\alpha)} = E^{\c_j(I_\alpha)-1}N &\text{ if } \c_j(I_\alpha) \geq 1 \text{ and } \d_{n_\alpha-j}(I_\alpha) = 0\,,\\
	&w_j^{(\alpha)} = WS^{\d_{n_\alpha-j}(I_\alpha)-1} &\text{ if } \c_j(I_\alpha) =0 \text{ and } \d_{n_\alpha-j}(I_\alpha) \geq 1\,,\\
	&w_j^{(\alpha)} = V &\text{ if } \c_j(I_\alpha) =0 \text{ and } \d_{n_\alpha-j}(I_\alpha) = 0\,.
	\end{align*}
	We can compute contact and descent vectors of $I$ using Proposition \ref{statsint}:
	\begin{align*}
	\Cont(I) &= (\c_0(I_2)+i,\c_0(I_1)+1-i,\c_1(I_1),...,\c_{n_1}(I_1),\c_1(I_2),...,\c_{n_2}(I_2))\,,\\
	\Drun(I) &= (\d_0(I_2),...,\d_{n_2-1}(I_2),\d_0(I_1)+1,\d_1(I_1),...,\d_{n_1}(I_1),\d_{n_2}(I_2))\,.
	\end{align*}
	We define the words $w_j$ corresponding to each component of those vectors:
	\begin{align*}
	&w_j = E^{\c_j(I)-1}N &\text{ if } \c_j(I) \geq 1 \text{ and } \d_{n-j}(I) = 0\,,\\
	&w_j = WS^{\d_{n-j}(I)-1} &\text{ if } \c_j(I) =0 \text{ and } \d_{n-j}(I) \geq 1\,,\\
	&w_j = V &\text{ if } \c_j(I) =0 \text{ and } \d_{n-j}(I) = 0\,.
	\end{align*}
	We directly have that $w_0 = E^{\c_0(I_2)+i-1}N$, $w_{j+1} = w_j^{(1)}$ for $1 \leq j \leq n_1-1$ and $w_{n_1+1+j} = w_j^{(2)}$ for $1 \leq j \leq n_2$. The values of $w_1$ and $w_{n_1+1}$ depend on the emptiness of $I_1$ and on $i$.
	\begin{itemize}
		\item If $I_1$ is empty, then $n_1+1 = 1$, $i=1$ and $w_1 = W$, and we get\\
		$E w_0...w_n S = E^{\c_0(I_2)+1}NW w_1^{(2)}...w_{n_2}^{(2)} S = ENWS \odot \Phi(I_2) = \Phi(\oplusI([\bullet,\bullet]^{\bullet 1},I_2)) = \Phi(I)\,.$
		\item If $I_1$ is not empty, then $n_1+1 > 1$ and $w_{n_1+1} = w_{n_1}^{(1)}S$. We consider two cases for the value of $i$:
		\begin{itemize}
			\item If $i \leq \c_0(I_1)$, then $w_1 = E^{c_0(I_1)-i}N$ and we get 
			\begin{align*}Ew_0...w_nS &= E^{\c_0(I_2)+i}NE^{\c_0(I_1)-i}Nw_1^{(1)}...w_{n_1}^{(1)}Sw_1^{(2)}...w_{n_2}^{(2)}S\\
			&= E^iNE^{\c_0(I_1)-i}Nw_1^{(1)}...w_{n_1}^{(1)}S^2 \odot \Phi(I_2)\\
			&= \maltese_i(\Phi(I_1)) \odot \Phi(I_2)\\
			&= \Phi(I)\,.
			\end{align*}
			\item If $i = \c_0(I_1) +1$, then $w_1 = V$ and we get 
			\begin{align*}Ew_0...w_nS &= E^{\c_0(I_2)+\c_0(I_1)+1}NVw_1^{(1)}...w_{n_1}^{(1)}Sw_1^{(2)}...w_{n_2}^{(2)}S\\
			&= E^{\c_0(I_1)+1}NVw_1^{(1)}...w_{n_1}^{(1)}S^2 \odot \Phi(I_2)\\
			&= \maltese_{\c_0(I_1)+1}(\Phi(I_1)) \odot \Phi(I_2)\\
			&= \Phi(I)\,.
			\end{align*}
		\end{itemize}
	\end{itemize}
\end{proof}

With this direct description, we are able to prove that $\Phi$ is compatible with conjugations and it is also easier to prove statistics correspondences under the bijection:

\begin{prop}
	For every Tamari interval $I$, we have $\ConjF(\Phi(I)) = \Phi(\ConjI(I))$, and the following relations hold:
	\begin{align*}
	&|\Phi(I)|_E = \Val(P)+1\,,&|\Phi(I)|_N = \DR(P)+1\,,\\
	&|\Phi(I)|_W = \Val(Q)+1\,,&|\Phi(I)|_S = \DR(Q)+1\,.		
	\end{align*}
\end{prop}

\begin{proof}
	Conjugation of Tamari intervals exchanges contact and descent vectors, so its action on corresponding extended fighting fish consists in reverting the word and making the exchanges of steps $E \leftrightarrow S$ and $N \leftrightarrow W$, hence the relation about conjugations.\\
	Let $I = [P,Q]$ be a Tamari interval of size $n$. We recall that $\Cont(I) = \Cont(P)$ and $\Drun(I)=\Drun(Q)$, and that both vectors have a component sum equal to $n$. Proposition~\ref{valdr} allows us to compute:
	\begin{align*}
	|\Phi(I)|_E = 1 + \sum_{\substack{0 \leq i \leq n \\ \c_i(P)\geq 1}} (\c_i(P)-1) &= 1 + \sum_{0 \leq i \leq n} (\c_i(P)-1) - \sum_{\substack{0 \leq i \leq n \\ \c_i(P)=0}} (\c_i(P)-1)\\
	&= \Big|\{0\leq i\leq n, \c_i(P) = 0\}\Big|\\
	&= \Val(P) +1\,,
	\end{align*}
	\[|\Phi(I)|_N = \sum_{\substack{0 \leq i \leq n \\ \c_{i}(P)\geq 1}} 1 = \Big|\{0\leq i\leq n, \c_i(P) \geq 1\}\Big| = \DR(P)+1\,.\]
	
	We then get the two additional equalities using the conjugations relation:
	\[|\Phi(I)|_W = |\ConjF(\Phi(I))|_N = |\Phi(\ConjI(I))|_N = \DR(\ConjD(Q))+1 = \Val(Q)+1\,,\]
	\[|\Phi(I)|_S = |\ConjF(\Phi(I))|_E = |\Phi(\ConjI(I))|_E = \Val(\ConjD(Q))+1 = \DR(Q)+1\,.\]
\end{proof}

\subsection{A formula for the Tamari distance}\label{subsec4.2}

Another nice feature of the bijection $\Phi$ is that we get an interpretation of the Tamari distance of an interval in terms of the area of the corresponding extended fighting fish. The area can be computed with the longitude and latitude functions, and so we get a natural way to express the Tamari distance of an interval in terms of the components of contact and descent vectors.

\begin{thm}\label{formula}
	Let $I$ be a Tamari interval of size $n$. Then its Tamari distance is given by:
	\[\L(I) = \sum_{0 \leq i < j \leq n} (\c_i(I)-1)(1-\d_{n-j}(I))\,.\]
\end{thm}

\begin{proof}
	The area of $\FF{F} \in \mathcal{EF}$ is
	$\Area(\FF{F}) = \sum_{i=1}^{m} (\lat(\FF{F}^{\leq i})-\lat(\FF{F}^{\leq i-1}))\Long(\FF{F}^{\leq i})$, where $m$ is the total number of steps composing $\FF{F}$. Then, if we write $\Phi(I) = Ew_0...w_nS$ as in Theorem 4, with $m_j$ being the total length of $w_j$, we can expand $\Area(\Phi(I))$ as
	$$\Area(\Phi(I)) = \sum_{j=0}^n \sum_{i=1}^{m_j} (\lat(w_j^{\leq i})-\lat(w_j^{\leq i-1}))(1+\Long(w_0)+...+\Long(w_{j-1})+\Long(w_j^{\leq i}))$$
	We have $\Long(w_k)=\c_k(I)-1$ and we need to distinguish between three cases for $w_j$:
	\begin{itemize}
		\item If $w_j = E^{\c_j(I)-1}N$, then $m_j=\c_j(I)$, $\lat(w_j^{\leq i})-\lat(w_j^{\leq i-1})$ is 0 for $i < \c_j(I)$ and 1 for $i=\c_j(I)$, and $\Long(w_j^{\leq \c_j(I)}) = \Long(w_j)= \c_j(I)-1$, so that we get the inner sum of the expression above equal to $(1-\d_{n-j}(I))(\sum_{k=0}^{j-1} (\c_k(I)-1) + \c_j(I))$.
		\item If $w_j = WS^{\d_j(I)-1}$, then $m_j=\d_{n-j}(I)$, $\lat(w_j^{\leq i})-\lat(w_j^{\leq i-1})$ is 0 for $i =1$ and $-1$ for $i>1$, and $\Long(w_j^{\leq i}) = -1=\c_j(I)-1$ and we still get the inner sum equal to $(1-\d_{n-j}(I))(\sum_{k=0}^{j-1} (\c_k(I)-1) + \c_j(I))$.
		\item If $w_j = V$, then $m_j=1$, $\lat(w_j^{\leq 1})-\lat(w_j^{\leq 0}) = 1$ and $\Long(w_j^{\leq 1}) = -1$ and the inner sum is once again $(1-\d_{n-j}(I))(\sum_{k=0}^{j-1} (\c_k(I)-1) + \c_j(I))$.
	\end{itemize}
	We now use the fact that $\c_j(I) \neq 0$ implies $\d_{n-j}(I) =0$ (Proposition \ref{compzero}) to get	
	\begin{align*}
	\Area(\Phi(I)) &= \sum_{j=0}^n (1-\d_{n-j}(I))(\sum_{k=0}^{j-1} (\c_k(I)-1) + \c_j(I))\\
	&= \sum_{0 \leq i < j \leq n} (\c_i(I)-1)(1-\d_{n-j}(I)) + \sum_{j=0}^{n} (1-\d_{n-j}(I))\c_j(I)\\
	&= \sum_{0 \leq i < j \leq n} (\c_i(I)-1)(1-\d_{n-j}(I)) + \sum_{j=0}^n \c_j(I)\\
	&= \sum_{0 \leq i < j \leq n} (\c_i(I)-1)(1-\d_{n-j}(I)) + n\,.
	\end{align*}
	We can now conclude with the equality $\L(I) = \Area(\Phi(I))-n$.
\end{proof}

We want to point out that we can prove Theorem \ref{formula} in a shorter way using the decomposition of Tamari intervals, because it provide us nice expressions of contact and descent vectors. Nonetheless, such a proof would give no insight on how the formula could have been guessed. This is why we preferred to present a constructive proof using the inherent structure of extended fighting fish, giving an interesting application of the bijection $\Phi$.\\

The distance preserving property of the conjugation of Tamari intervals can be seen easily with this formula. It is also interesting to note that $\sum_{0 \leq i < j \leq n} (\c_i(P)-1)(1-\d_{n-j}(Q))$ can be defined even if $P$ and $Q$ are not comparable, and that it may define a statistic generalizing the Tamari distance on pairs of Dyck paths of the same size. Note however that it may not be the exact sum to consider: when $P$ and $Q$ are comparable, $\sum_{i=0}^n (\c_i(P)-1)(1-\d_{n-i}(Q)) = n-1$, so we may consider the other Tamari distance generalizing statistic $\sum_{0 \leq i \leq j \leq n} (\c_i(P)-1)(1-\d_{n-j}(Q)) - (n-1)$, which does not agree in general with the sum in Theorem \ref{formula} when $P$ and $Q$ are not comparable. Several other identities over contact and descent vectors are satisfied when $P$ and $Q$ are comparable, so there are many other candidates for a statistic generalizing the distance of Tamari intervals.

\subsection{The average area of extended fighting fish}\label{subsec4.3}

\subsubsection{An equation for extended fighting fish}\label{subsubsec4.3.1}

Let us denote by $H(u,q)\equiv H(t;u,q,y)=\sum_{F \in \mathcal{EF}}t^{n(F)}u^{j(F)}q^{a(F)}y^{v(F)}$, where $n(F)=\size(F)$, $j(F)=\ljaw(F)$, $a(F)=\Area(F)$ and $v(F)=|F|_V$, the generating function of extended fighting fish according to their size, area, jaw length and number of vertical steps. We are mostly interested in the size and area but we need to keep track of the jaw length, which is a so-called catalytic parameter, playing an essential role in the description of the decomposition of Theorem \ref{decfish}. The variable $y$ allows us to deal simultaneously with the generating functions of extended fighting fish (for $y=1$) and that of fighting fish (for $y=0$) by keeping track of the number of vertical steps. Using the standard symbolic method~\cite[Chapter 1]{FS}, the decomposition of extended fighting fish can be translated into the following equation for the generating function:

\begin{align}\label{eqHuq}
H\left(u,q \right)=\,tuq \left( 1+H \left(u,q \right)  \right) +&\frac {tuq \left( H \left(1,q \right) -H \left( uq,q \right)  \right)  \left( 1+H \left(u,q \right)  \right) }{1-uq}\notag\\
&\hspace{123pt} +tuqyH \left(uq,q \right)  \left( 1+H \left(u,q \right)  \right)\,.
\end{align}

In particular, in the previous equation, the term $tuq \left( 1+H \left(u,q \right)  \right)$ corresponds to the case $F=\oplusF(\varepsilon,\FF{F}_2) = \maltese_i(\varepsilon) \odot \FF{F}_2$ of Theorem~\ref{decfish}, the term $\frac {tuq \left( H \left(1,q \right) -H \left( uq,q \right)  \right)  \left( 1+H \left(u,q \right)  \right) }{1-uq}$ corresponds to the case $F=\oplusF(\FF{F}_1^{\bullet i},\FF{F}_2) = \maltese_i(\FF{F}_1) \odot \FF{F}_2$  with  $1\leq i\leq\ljaw(F_1)$ and the term $tuqyH \left(uq,q \right)  \left( 1+H \left(u,q \right)  \right)$ corresponds to the case $F=\oplusF(\FF{F}_1^{\bullet i},\FF{F}_2) = \maltese_i(\FF{F}_1) \odot \FF{F}_2$ with $i=\ljaw(F_1)+1$.

In order to first obtain $H(1,1)$, the generating function of extended fighting fish according to their size, we set $q=1$ to obtain the following equation, which is polynomial but not algebraic because it still involves  two unknown series $H(1,1)$ and $H(1,u)$:
$$
-(1-u)H\left(u,1 \right)+tu \left(H \left(1,1 \right) -H \left( u,1 \right) +(1-u)\left(1+yH \left(u,1 \right)\right)\right)  \left( 1+H \left(u,1 \right)  \right)=0
$$
According to \cite{MBMA}, we can solve this equation by applying the extended kernel method in order to obtain an algebraic equation for $H(1,1)$. Indeed, upon deriving the previous equation with respect to $u$ we obtain the following:
\begin{align}\label{bigeq}
& \bigg(
-(1-u)+tu \Big(H \left(1,1 \right) -H \left( u,1 \right) +(1-u)\big(1+yH \left(u,1 \right)\big)\Big)\notag\\
&\qquad\qquad\qquad\qquad\qquad\qquad\qquad\qquad\qquad\qquad
+tu \big(-1 +(1-u)y\big)  \big( 1+H \left(u,1 \right)  \big)\bigg)
\frac{\partial H}{\partial u}  \left(u,1 \right)\notag\\
& + H\left(u,1 \right)+t \Big(H \left(1,1 \right) -H \left( u,1 \right) +(1-u)\big(1+yH \left(u,1 \right)\big)\Big)  \big( 1+H \left(u,1 \right)  \big)\notag\\
&\qquad\qquad\qquad\qquad\qquad\qquad\qquad\qquad\qquad\qquad+tu \Big(-\big(1+yH \left(u,1 \right)\big)\Big)  \big( 1+H \left(u,1 \right)  \big) = 0\,.
\end{align}

If we find a power series $U\equiv U(t)$ that can be substituted in the equation \eqref{bigeq} to cancel the coefficient of $\frac{\partial H}{\partial u} (u,1)$, then the formal power series $U$, $H(U,1)$ and $H(1,1)$ would also cancel the rest of the equation, and together with the master equation~\eqref{eqHuq}, we would have a system of three polynomials for these three unknown series.

Now observe that the equation for $U$ cancelling the coefficient of $\partial H/\partial u$ can be rewritten as:
\begin{align*}
& U=1-tU \left(H \left(1,1 \right) -H \left( U,1 \right) +(1-U)\left(1+yH \left(U,1 \right)\right)+tU \left(-1 +(1-U)y\right)  \left( 1+H \left(U,1 \right)  \right)\right)\,.
\end{align*}
from which, upon observing that $H(u,q)$ is a formal power series in $Q[u,q][[t]]$, one can infer by recurrence that there indeed exists a unique such formal power series $U$.

By standard elimination techniques on the three equations, and upon setting $U=1+V$, we obtain
$$
\left\{
\begin{array}{cl}
V&=t\,(1+Vy)(1+V)^3\,,\\
H(1,1)&=V-V^2-V^3y\,.
\end{array}
\right.
$$

In particular, for $y=0$ we recover the known parametrization of the generating function of fighting fish or synchronized Tamari interval, and, for $y=1$ a parametrization for the generating function of extended fighting fish which corresponds to the standard parametrization for the generating function of Tamari intervals.

\subsubsection{The total area generating function}\label{subsubsec4.3.2}
The generating function of the total area of parametrized extended fighting fish is the series A(t) given by
$$
A\equiv A(t)=\frac{\partial H}{\partial q}(1,1)=\left(\frac{\partial}{\partial q}H(u,q)\right)\mid_{u=1,q=1}\,.
$$

In order to obtain it, we can differentiate the master equation with respect to $q$ and set $q=1$. Upon setting also $u=U$  a first simplification occurs because, in view of the chain rule for derivation, the coefficient of $\frac{\partial H}{\partial q}(u,1)$ in the derivative of the master equation \eqref{eqHuq} with respect to $q$ is the same, after putting $q=1$, as the coefficient of $\frac{\partial H}{\partial u}(u,1)$ in the derivative of the master equation with respect to $u$, which is canceled by $u=U$.  The remaining terms can be further simplified upon using $U=1+V$, $H(1,1)=V-V^2-v^3y$ and $H(U,1)=V$, and the resulting equation is:

\begin{align}\label{eqDerivate}
\frac{\partial H}{\partial q} (1,1)=(1+V)(1+yV)\cdot \frac{\partial H}{\partial u} (1+V,1)\,.
\end{align}

In order to compute $\frac{\partial H}{\partial u}(1+V,1)$ we restart from the derivative of the master equation with respect to $u$ at $q=1$ and derive a second time with respect to $u$: the coefficient of $\frac{\partial^2H}{\partial u^2}(u,1)$ in this second derivative of the master equation is again precisely the same as the coefficient of $\frac{\partial H}{\partial u}(u,1)$ in the first derivative of the master equation with respect to $u$ at $q=1$, which is canceled by $u=U=1+V$. The resulting equation is a quadratic equation for $\frac{\partial H}{\partial u} (1+V,1)$:
\begin{align*}
V(1+yV)\cdot\left(\frac{\partial H}{\partial u} (1+V,1)\right)^2+(V^2y-1)\cdot\frac{\partial H}{\partial u} (1+V,1)+V(1+yV)=0\,.
\end{align*}

In turn, using Equation~\eqref{eqDerivate}, this yields a quadratic equation for the total area generating function $A$:
\[
V A^2+(1+V)(V^2y-1)A+V(1+V)^2(1+yV)^2=0\,.
\]

In particular, for $y=0$, we recover the equation for the area generating function of fighting fish
\[
V A^2+(1+V)A+V(1+V)^2 = 0
\]
or
\[
A=\frac{1+V}{2V}\left(1-\sqrt{1-4V^2}\right)\,,
\]
while for $y=1$, we have
\[
V A^2+(1+V)^2(V-1)A+V(1+V)^4 = 0
\]
or
\[
A=\frac{1+V}{2V}\left(1-V^2-\sqrt{(1+V)^3(1-3V)}\right)\,.
\]

\subsubsection{Asymptotics}\label{subsubsec4.3.3}
The asymptotic behavior for fighting fish ($y=0$) is already given in~\cite{ff}. We give here 
the asymptotic behavior of the number of extended fighting fish and total area of extended fighting fish. In particular we obtain the following theorem:
\begin{thm}
  The average area of uniformly random extended fighting fish of size $n$ has the following asymptotic behavior when $n$ grows to infinity:
  \[
    \frac{[t^n]A(t)}{[t^n]H(t,1,1)}\sim
    \frac{2^{1/4}3^{3/4}\sqrt{\pi}}{2\Gamma({\frac3{4}})} n^{5/4}\,.
  \]
\end{thm}

This behavior in $n^{5/4}$ is the same as for fighting fish, and then belongs to a different universality class from the one containing all classical models of polyominoes for which the area grows like $n^{3/2}$. 

\begin{proof}
First following~\cite{FS}, Theorem VII.2, we obtain the singular expansion of $V$ around it unique dominant singularity $\rho=27/256$:
\[
V(t)=\frac13-\frac{2\sqrt 6}9\sqrt{1-t/\rho}+\frac{10}{27}(1-t/\rho)+O((1-t/\rho)^{3/2})\,.
\]

Using this expansion of $V$ we obtain the expansion of $H(1,1)$, which has a dominant singular term of order $3/2$ as expected for solutions of equations with one catalytic variable~\cite{DNY}
\[
H(t,1,1)=\frac{5}{27}-\frac{16}{27}(1-t/\rho)+\frac{32\sqrt{6}}{81}(1-t/\rho)^{3/2}+O((1-t/\rho)^2)\,,
\]
and that of $A(t)$, which displays instead the fish-area behavior with a dominant singular term of order $1/4$ (see \cite{ff}):
\[
A(t,1,1)=\frac{16}{9}-\frac{16\sqrt{2}\,6^{1/4}}{9}(1-t/\rho)^{1/4}+O(\sqrt{1-t/\rho})\,.
\]

Using the transfert theorems of~\cite{FS}[Theorem VI.2], we have
\[
[t^n]H(t,1,1)\sim \frac3{4\sqrt{\pi}}\frac{32\sqrt6}{81}\rho^{-n}\,n^{-5/2} \textrm{ and }[t^n]A(t)\sim\frac1{4\Gamma(\frac3{4})}\frac{16\sqrt{2}\,6^{1/4}}{9}\rho^{-n}n^{-5/4}\,,
\]
so that the average area of extended fighting fish is
\[
\frac{[t^n]A(t)}{[t^n]H(t,1,1)}\sim \frac{2^{1/4}3^{3/4}\sqrt{\pi}}{2\Gamma({\frac3{4}})} n^{5/4}\,.
\]
\end{proof}

\section{Final comments}\label{sec5}

Our bijection $\Phi$ gives a nice way to represent intervals of the Tamari lattice as paths in the quarter plane or branching surfaces. It is natural in the sense that statistics and structure of both objects are transferred by $\Phi$. Also, we want to note that the decomposition we presented here for both extended fighting fish and Tamari intervals can also be stated for bridgeless planar maps, and that the subdecomposition of fighting fish and synchronized intervals already appeared for nonseparable planar maps~\cite{Brown, fangnonsep}, left ternary trees~\cite{DDP} and two-stack sortable permutations~\cite{goulden}. From all these isomorphic decompositions can be derived recursive bijections with our (extended) fighting fish. It might be interesting to find which statistic corresponds to the area of fighting fish for the different mentioned structures.\\

In~\cite{bergeron2011higher}, the authors define an order on the set of $m$-Dyck paths ($m \in \mathbb{N}^*$) that gives rise to a lattice called the $m$-Tamari lattice. Indeed, the case $m=1$ corresponds to the classical Tamari lattice that we defined in Section~\ref{sec3}. Intervals in the $m$-Tamari lattice of order n are still counted~\cite{BMFPR11} by a beautiful and simple number: $\frac{m+1}{n(mn+1)}\binom{(m+1)^2n+m}{n-1}$. However, no bijective proof for this formula is known, and finding a model of fish for $m$-Tamari intervals could be a first step towards this direction. Let us also point out that the formula of Theorem~\ref{formula} for the Tamari distance can be extended to a formula for the distance of $m$-Tamari intervals since the $m$-Tamari lattice of order $n$ can be seen as an upper ideal of the $1$-Tamari lattice of order $mn$.

\section*{Acknowledgements}

We would like to thank the anonymous referees for the comments that helped us improve the precision and the clarity of the article.

\bibliographystyle{elsarticle-harv}
\bibliography{EJC-v5}

\end{document}